\documentclass[11pt]{article} \topmargin=-0.5cm
\usepackage{amssymb}
\def\citet{\cite}

   \csname @addtoreset\endcsname{equation}{section}
    \csname @addtoreset\endcsname{figure}{section}
    \csname @addtoreset\endcsname{table}{section}

\newcounter{thanksnum}
\def\thanksnumber#1
{\setcounter{thanksnum}{\value{footnote}}\setcounter{footnote}{#1}%
                     \addtocounter{footnote}{-1}\footnotemark
                     \setcounter{footnote}{\value{thanksnum}}}

\newtheorem{theorem}{Theorem}[section]
\newtheorem{lemma}{Lemma}[section]

\newtheorem{proposition}{Proposition}[section]

\newtheorem{condition}{Condition}[section]
\newtheorem{definition}{Definition}[section]
\newtheorem{remark}{Remark}[section]

\def\e{\varepsilon}

\def\defi{\stackrel{{\scriptscriptstyle \Delta}}{=}}

\def\a{\alpha}
\def\d{\delta}
\def\o{\omega}
\def\O{\Omega}

\def\Y{{\cal Y}}
\def\F{{\cal F}}
\def\w{\widehat}
\def\Ind{{\mathbb{I}}}

\def\esssup{\mathop{\rm ess\, sup}}
\def\essinf{\mathop{\rm ess\, inf}}
\def\const{{\rm const\,}}

\def\R{{\bf R}}
\def\E{{\bf E}}
\def\P{{\bf P}}

\def\Z{{\cal Z}}
\def\PP{{\cal P}}

\def\H{{\cal H}}

\def\L{L}

\def\b{\beta}
\def\s{\delta}
\def\g{\gamma}

\def\W{{\cal W}^*}

\def\X{{\cal X}}
\def\t{\theta}
\def\oo{\bar}
\def\s{\sigma}
\def\D{{\Delta}}
\def\p{\partial}

\def\U{{\cal U}}
\def\V{{\cal V}}
\def\A{{\cal A}}
\def\M{{\cal M}}
\def\B{{\cal B}}
\def\L{{\cal L}}

\newcommand{\be}{\begin{equation}}
\newcommand{\ee}{\end{equation}}
\newcommand{\bd}{\begin{displaymath}}
\newcommand{\ed}{\end{displaymath}}
\newcommand{\ba}{\begin{array}{ll}}
\newcommand{\ea}{\end{array}}
\newcommand{\baa}{\begin{eqnarray}}
\newcommand{\eaa}{\end{eqnarray}}
\newcommand{\baaa}{\begin{eqnarray*}}
\newcommand{\eaaa}{\end{eqnarray*}}   \font\sm=cmr10

\def\W{{\cal W}}
\def\H{{\cal H}}

\def\H{{\cal H}}

\def\CC{{\cal C}}
\def\CC{{ \texttt{C}}}
\def\WW{\varrho}
\def\A{A}
\date{ Submitted November 1, 2016. Revised June  21, 2017}
\title{On recovering solutions for SPDEs from their averages}
\author{
Nikolai Dokuchaev\\
 {\sm Department of Mathematics \& Statistics, Curtin
University,}\\ {\sm  GPO Box U1987, Perth, 6845 Western Australia} }
\begin{document}
\maketitle
\begin{abstract}  We study linear stochastic partial differential equations of
parabolic type. We consider a new boundary value problem where a Cauchy condition is replaced by
 a prescribed average of the solution either over time and probabilistic space  for forward SPDEs and over time for backward SPDEs.
 Well-posedness, existence, uniqueness,
 and  a regularity of the solution for this new problem are obtained.
 In particular, this can be considered as a possibility to recover a solution of a forward SPDE
 in a setting where its values at the initial time are  unknown, and where  the average
 of the solution over time and probability space is observable, as well as the  input processes.
\end{abstract}
\section{Introduction}
The paper studies boundary value problems for
stochastic partial differential equations of the second order.
These equations have many applications and were widely studied; see e.g.
 \cite{Alos,Ba,DaPT,D92,D95,D05,D10,D12,DT, A1,kl,kr,Ma,Mas,Mat,Moh,Par,Roz,Wa,Yo,Zh}.
\index{Al\'os et al (1999), Bally {\it et al} (1994),  Da Prato and Tubaro
(1996), Dokuchaev (1992,2005,2010,2012), Du and Tang (2012), Gy\"ongy (1998), Krylov (1999),  Ma and
Yong (1997), Maslowski (1995), Pardoux (1993),
 Rozovskii (1990), Walsh (1986),  Yong and Zhou
 (1999), Zhou (1992).}
Forward parabolic SPDEs  are usually considered with a Cauchy
condition at initial time, and backward parabolic SPDEs are usually
considered with a Cauchy condition at terminal time. However,
there are also results for SPDEs with boundary conditions that mix
the solution at different times that may  include initial time and
terminal time. This category includes stationary type solutions for
forward SPDEs; see, e.g.,  \cite{CK, CM,CMG,Duan,liu,Mat,Moh,Sin}, \index{Caraballo {\em et al } (2004),
 Chojnowska-Michalik (1987), Chojnowska-Michalik and Goldys
 (1995), Duan {\em et al} (2003), Mattingly (1999),
Mohammed  {\em et al} (2008), Sinai (1996),} and the references
therein.  Related results were obtained for periodic solutions of
SPDEs in \cite{CM90,FZ,kl}.\index{
 (Chojnowska-Michalik (1990), Feng and Zhao (2012), Kl\"unger  (2001)).}
 Some results for parabolic equations and stochastic PDEs
with non-local conditions replacing the Cauchy condition  were obtained  in
\citet{D94,D04,D08,D11,D15}. \index{Dokuchaev (2008,2011,2015).}

The present paper addresses these and related problems again.
We
consider forward and backward SPDEs with the Dirichlet  condition at the boundary of the
state domain; the equations are of  a parabolic type. For forward SPDEs, a  Cauchy condition at initial time  is replaced by
 a condition requiring a prescribed average of the solution over time and probabilistic space. For backward SPDEs, a  Cauchy condition at
 terminal time  is replaced by
 a condition requiring a prescribed average of the solution over time.
 This is a novel setting for SPDEs; for deterministic parabolic equations, a related result was obtained in \citet{D16}.    We
obtained sufficient conditions for existence and regularity of the
solutions in $L_2$-setting  (Theorems \ref{FWTh1} and \ref{ThB} below).
 This result can be interpreted as a possibility to recover a parabolic diffusion
from its time-average when the values at the initial time are  unknown.
\section{The problem setting and definitions}
Assume that  we
are given a standard  complete probability space $(\O,\F,\P)$ and a
right-continuous filtration $\F_t$ of complete $\s$-algebras of
events, $t\ge 0$, such that  $\F_0$ is the $\P$-augmentation of the set
$\{\emptyset,\O\}$. We are given also a $N$-dimensional Wiener process
$w(t)$ with independent components;  it is a Wiener process with
respect to $\F_t$.
\par
Let  $D\subset\R^n$ be an open bounded domain with
a $C^2$-smooth boundary $\p D$. Let $T>0$ be given, and let $Q\defi
D\times [0,T]$. \par
 We consider the following boundary value problem in $Q$ for forward SPDEs
\baa
\label{FWparab1} &&d_tu=(\A u+ \varphi)\,dt +\sum_{i=1}^N
[B_iu+h_i]\,dw_i(t), \quad t\ge 0,
\\\label{FWparab10}
&& u(x,t,\o)\,|_{x\in \p D}=0,
\\ &&\E\left(\kappa  u(x,T)+\int_0^T\WW(t)u(x,t)dt\right)=\mu(x).\label{FWparab2}
\eaa

Here  $\mu$, $\varphi$, and $h_i$ are given inputs, $u$ is a sought out solution.

In addition, we will study the following boundary
value problem in $Q$ for backward SPDEs
\begin{eqnarray} 
\label{parab1} &&d_tu+(\A u+ \varphi)\,dt +\sum_{i=1}^N
B_i\chi_idt=\sum_{i=1}^N\chi_i dw_i(t), \quad t\ge 0,
\\\label{parab10}
&& u(x,t,\o)\,|_{x\in \p D}=0
\\ && \kappa u(x,0)+\int_0^T\WW(t)u(x,t,\o)dt=\psi(x,\o)\quad \hbox{a.s. for}\quad x\in D.\
\label{parab2}
\end{eqnarray}

Here  $\mu$ and  $\varphi$ are given inputs, and the set of functions $(u,\chi_1,...,\chi_N)$ is a sought out solution.

Here $u=u(x,t,\o)$, $\psi=\psi(x,\o)$, $\varphi=\varphi(x,t,\o)$, $h_i=h_i(x,t,\o)$,
$\chi_i=\chi_i(x,t,\o)$,
 $(x,t)\in Q$,   $\o\in\O$.\par

In (\ref{FWparab2}), $\kappa\in \R$,  $\WW(t)$ is a measurable and bounded non-random function.

In these SPDEs, $\A$ and $B$ are differential operators defined as
\baaa
&&\A v\defi \sum_{i=1}^n  \frac{\p }{\p x_i}\left(\sum_{j=1}^n a_{ij}(x)\frac{\p v}{\p x_j}(x)\right)+a_0(x)v(x),
\\
\label{B}
&&B_iv\defi\frac{dv}{dx}\,(x)\,\oo\beta_i(x,t,\o) +
\beta_i(x,t,\o)\,v(x),\quad i=1,\ldots ,N, \eaaa
where  functions  $\b_j(x,t,\o):
\R^n\times[0,T]\times\O\to\R^n$, $\oo\b_i(x,t,\o):$
$\R^n\times[0,T]\times\O\to\R$, and $\varphi (x,t,\o): \R^n\times
[0,T]\times\O\to\R$ are progressively measurable with respect to
$\F_t$ for all $x\in\R^n$. The function $\psi(x,\o):
\R^n\times\O\to\R$ is $\F_T$-measurable for all $x\in\R^n$.
\subsection*{Conditions for the coefficients}
 To proceed further, we assume that Conditions
\ref{cond1}-\ref{condK} remain in force throughout this paper.
 \begin{condition}\label{cond1}
 The functions $a_{ij}(x): D\to \R$ and $a_0(x): D\to
 \R$ are continuous and bounded, and there exist
 continuous bounded derivatives  $\p a_0(x)/\p x_i$, $\p a_{ij}(x)/\p x_i$, $i,j=1,...,n$.  In
addition, we assume that the matrix $a=\{a_{ij}\}$ is symmetric.The functions
$\oo\b_i(x,t,\o)$ and $\b_i(x,t,\o)$ are bounded and
differentiable in $x$ for a.e. $t,\o$, and the corresponding
derivatives are bounded.
\end{condition}
It follows from this condition that there exist modifications of
$\b_i$ such that the functions $\b_i(x,t,\o)$ are continuous in $x$
for a.e. $t,\o$. We assume that $\b_i$ are such functions.
\begin{condition}\label{condB}
 $\oo\b_i(x,t,\o)=0$ for $x\in \p D$,
$i=1,...,N$.
\end{condition}
\begin{condition}\label{condA} [Superparabolicity condition \citet{Roz}] There exists a constant $\d>0$ such that \be
 \label{Main1} y^\top  a
(x)\,y-\frac{1}{2}\sum_{i=1}^N |y^\top\oo\b_i(x,t,\o)|^2 \ge
\d|y|^2 \quad\forall\, y\in \R^n,\ (x,t)\in  D\times [0,T],\
\o\in\O. \ee
\end{condition}
\par
 If $\kappa\neq 0$ and $\WW\equiv 0$, then the boundary value  problems above
 are ill-posed, with ill-posed  Cauchy conditions $u(x,T)=\mu(x)$ or  $u(x,0)=\psi(x)$
 respectively.
This case is excluded from consideration by imposing the following restrictions.
\begin{condition}\label{condK}
\begin{enumerate}
\item $\WW(t)\ge 0$ a.e. and $\kappa\ge 0$.
\item For problem (\ref{FWparab1})-(\ref{FWparab2}), we assume  that there exists $T_1\in(0,T]$ such that \baaa
\essinf_{t\in[0,T_1]}\WW(t)>0.\eaaa
\item
For problem (\ref{parab1})-(\ref{parab2}), we assume   that the function $\WW(t)$ is continuous at $t=T$
 and that
there exists $T_1\in[0,T)$ such that \baaa
\essinf_{t\in[T_1,T]}\WW(t)>0.\eaaa
\end{enumerate}
\end{condition}
\par
We do not exclude an important special case where the function
 $\varphi$ is deterministic, and $h_i\equiv 0$, $B_i\equiv 0$ for all $i$. In this case, the boundary value problem is deterministic,
 and   $\chi_i\equiv 0$ ($\forall i)$ for backward equations.

\subsection*{Spaces and classes of functions} 
 We denote  by $|\cdot|$ the Euclidean norm in $\R^k$, and $\bar G$ denote
the closure of a region $G\subset\R^k$.
We denote by $\|\cdot\|_{ X}$ the norm in a linear normed space $X$,
and
 $(\cdot, \cdot )_{ X}$ denote  the scalar product in  a Hilbert space $
X$.
\par
Let us introduce some spaces of real valued functions.
\par
 Let $G\subset \R^d$ be an open
domain. For $q\ge 1$, we denote by $L_q(G)$ the usual Banach spaces
of classes  of equivalency of measurable by Lebesgue functions
$v:G\to \R$, with the norms
$\|v\|_{L_q(G)}=\left(\int_G|v(x)|^qdx\right)^{1/q}$.  For integers
$m\ge 0$, we denote by ${W_q^m}(G)$ the Sobolev  spaces of functions
that belong to $L_q(G)$ together with the distributional derivatives
up to the $m$th order, $q\ge 1$, with the norms
$\|v\|_{W_q^m(G)}=\left(\sum_{k:\, |k|\le m}\|D^{k}v\|_{L_q(G)}^q
\right)^{1/q}$. 
\def\D{{\cal D}}
Here $\D^k=\D_{k_1}\cdots \D_{k_d}$ is the partial
derivative of the order $|k|=\sum_{i=1}^d k_i$, where $k=(k_1,...,k_d)\in\mathbb{Z}^d$,
 $k_i\ge 0$, $\D_i=\p^{k_i}/dx_k^{k_i}$.

\par Let $H^0\defi L_2(D)$,
and let $H^1\defi \stackrel{\scriptscriptstyle 0}{W_2^1}(D)$ be the
closure in the ${W}_2^1(D)$-norm of the set of all smooth functions
$u:D\to\R$ such that  $u|_{\p D}\equiv 0$. Let $H^2=W^2_2(D)\cap
H^1$ be the space equipped with the norm of $W_2^2(D)$. The spaces
$H^k$  are Hilbert
spaces, and $H^k$ is a closed subspace of $W_2^k(D)$, $k=1,2$.
\par
 Let $H^{-1}$ be the dual space to $H^{1}$, with the
norm $\| \,\cdot\,\| _{H^{-1}}$ such that if $u \in H^{0}$ then $\|
u\|_{ H^{-1}}$ is the supremum of $(u,v)_{H^0}$ over all $v \in H^1$
such that $\| v\|_{H^1} \le 1 $. $H^{-1}$ is a Hilbert space.
\par
Let $C_0(\oo D)$ be the Banach space of all functions $u\in C(\oo
D)$  such that $u|_{\p D}\equiv 0$ equipped with the norm from
$C(\oo D)$.
\par We shall write $(u,v)_{H^0}$ for $u\in H^{-1}$
and $v\in H^1$, meaning the obvious extension of the bilinear form
from $u\in H^{0}$ and $v\in H^1$.
\par
We denote by $\oo\ell _{k}$ the Lebesgue measure in $\R^k$, and we
denote by $ \oo{\B}_{k}$ the $\sigma$-algebra of Lebesgue sets in
$\R^k$.
\par
We denote by $\oo{{\PP}}$  the completion (with respect to the
measure $\oo\ell_1\times\P$) of the $\s$-algebra of subsets of
$[0,T]\times\O$, generated by functions that are progressively
measurable with respect to $\F_t$.
\par
 We  introduce the spaces
 $$
{\cal C}_k\defi C\left([s,T]; H^k\right),\quad \W^{k}\defi
L^{2}\bigl([ 0,T ],\oo{\cal B}_1, \oo\ell_{1};  H^{k}\bigr), \quad  k=-1,0,1,2,
$$
and the spaces
$$
\V^k(s,T)\defi \W^{1}(s,T)\cap {\cal C}_{k-1},\quad k=1,2,
$$
with the  norm $ \| u\| _{\V} \defi \| u\| _{{\W}^k} +\|
u\| _{{\cal C}_{k-1}}. $

For $\t\in [0,T]$, we introduce a space $\U_\t$ of functions $\varphi\in\W^0$ such that $\varphi(\cdot,t)=\varphi(\cdot,\t)+\int_\t^t\w\varphi(\cdot,s)ds$ for $t\in[\t,T]$ for
some $\w u\in L_1([\t,T];H^0)$, with the norm
 \baaa\| \varphi\| _{\U_\t} \defi \| \varphi\| _{{\W}^0} +\|\varphi(\cdot,\t)\|_{H^0}+ \int_\t^T\|\w\varphi(\cdot,t)\|_{H^0}dt. \eaaa
In particular, $\varphi(\cdot,t)$ is continuous in $H^0$ in $t\in[T-\t,T]$.
 If $\t=T$ then $\U_\t=\W^0=L_2(D\times[0,T])$.

In addition, we introduce the spaces
\baaa
 &&X^{k}\defi L^{2}\bigl([0,T]\times\Omega,
{\oo{\PP }},\oo\ell_{1}\times\P;  H^{k}\bigr), \quad\\ &&Z^k_t
\defi L^2\bigl(\Omega,{\F}_t,\P; H^k\bigr),\\
&&\CC^{k}_Z\defi C\left([0,T]; Z^k_T\right), \qquad k=-1,0,1,2,
 \eaaa
and the spaces $$ Y^{k}\defi
X^{k}\!\cap \CC^{k-1}_Z, \quad k=1,2, $$ with the norm $ \|
u\| _{Y^k(s,T)}
\defi \| u\| _{{X}^k} +\| u\| _{\CC^{k-1}_Z}. $

We introduce Banach spaces $\Z_p^k=L_p(\O,\F,\P,H^k)$, $p\in[1,+\infty]$.

  The spaces $\W_k$, $X^k$, and $Z_t^k$,  are Hilbert spaces.

\par
\index{For a set $S$ and a normed space ${\rm X}$, we denote by
${\rm B}(S,{\bf X})$ the set of all bounded functions $x:S\to{\rm
X}$. For a set $S$ and a Banach space ${\rm X}$, we denote by ${\rm
B}(S,{\bf X})$ the Banach space of bounded functions $x:S\to{\rm X}$
equipped with the norm $\|u\|_{{\rm B}}=\sup_{s\in S}\|x(s)\|_{{\rm
X}}$. We introduce space $\Z=Z_T^0\cup B(D\times\O)$ equipped with
the norm $\|z\|_{\Z}= \|z\|_{Z_T^0}+\|z\|_{ B(D\times\O) }$.}

\begin{proposition} 
\label{propL} Let $\zeta\in X^0$, and
 let a sequence  $\{\zeta_k\}_{k=1}^{+\infty}\subset
L^{\infty}([0,T]\times\O, \ell_1\times\P;\,C(D))$ be such that all
$\zeta_k(\cdot,t,\o)$ are progressively measurable with respect to
$\F_t$, and  $\|\zeta-\zeta_k\|_{X^0}\to 0$ as $k\to +\infty$. Let
$t\in [0,T]$ and $j\in\{1,\ldots, N\}$ be given.
 Then the sequence of the
integrals $\int_0^t\zeta_k(x,s,\o)\,dw_j(s)$ converges in $Z_t^0$ as
$k\to\infty$, and its limit depends on $\zeta$, but does not depend
on $\{\zeta_k\}$.
\end{proposition}
\par
{\em Proof} follows from completeness of  $X^0$ and from the
equality
\begin{eqnarray*}
&&\E\int_0^t\|\zeta_{k}(\cdot,s,\o)-\zeta_m(\cdot,s,\o)\|_{H^0}^2\,ds
=\int_D\,dx\,\E\left[\int_0^t\big(\zeta_k(x,s,\o)-
\zeta_m(x,s,\o)\big)\,dw_j(s)\right]^2.
\end{eqnarray*}
\begin{definition} 
\rm For $\zeta\in X^0$, $t\in [0,T]$, $j\in\{1,\ldots, N\}$,  we
define $\int_0^t\zeta(x,s,\o)\,dw_j(s)$ as the limit  in $Z_t^0$ as
$k\to\infty$ of a sequence $\int_0^t\zeta_k(x,s,\o)\,dw_j(s)$, where
the sequence $\{\zeta_k\}$ is such  as in Proposition \ref{propL}.
\end{definition}

\par
Sometimes we shall omit $\o$.

\subsection*{The definition of solution for forward SPDEs}

\begin{definition} 
\label{FWdefsolltion} \rm Let $u\in Y^1$,  $\varphi\in X^{-1}$, and
$h_i\in X^0$. We say that equations
(\ref{FWparab1})-(\ref{FWparab10}) are satisfied if
\begin{eqnarray}
u(\cdot,t,\o)=u(\cdot,0,\o)&+& \int_0^t\big(\A u(\cdot,s,\o)+
\varphi(\cdot,s,\o)\big)\,ds  \nonumber
\\  \hphantom{xxxxxxxxx}&+&  \sum_{i=1}^N
\int_0^t[B_iu(\cdot,s,\o)+h_i(\cdot,s,\o)]\,dw_i(s)
\label{FWintur}
\end{eqnarray}
for all $t\in [0,T]$, and this equality is satisfied as an equality
in $Z_T^{-1}$.
\end{definition}
\par
Note that the condition on $\p D$ is satisfied in the  sense that
$u(\cdot,t,\o)\in H^1$ for a.e. \ $t,\o$. Further, $u\in Y^1$, and
the value of  $u(\cdot,t,\o)$ is uniquely defined in $Z_T^0$ given
$t$, by the definitions of the corresponding spaces. The integrals
with $dw_i$ in (\ref{FWintur}) are defined as elements of $Z_T^0$.
The integral with $ds$ in (\ref{FWintur}) is defined as an element
of $Z_T^{-1}$. In fact, Definition \ref{FWdefsolltion} requires for
(\ref{FWparab1}) that this integral must be equal  to an element of
$Z_T^{0}$ in the sense of equality in $Z_T^{-1}$.

\subsection*{The definition of solution for backward SPDEs}\begin{definition} 
\label{defsolltion} \rm Let $u\in Y^1$, $\chi_i\in X^0$,
$i=1,...,N$, and $\varphi\in X^{-1}$. We say that equations
(\ref{parab1})-(\ref{parab10}) are satisfied if \baa
u(\cdot,t,\o)=u(\cdot,T,\o)&+& \int_t^T\big(\A u(\cdot,s,\o)+
\varphi(\cdot,s,\o)\big)\,ds \ \nonumber
\\\hphantom{xxx}&+& \sum_{i=1}^N
\int_t^TB_i\chi_i(\cdot,s,\o)ds-\sum_{i=1}^N
\int_t^T\chi_i(\cdot,s)\,dw_i(s)
\label{intur} \eaa for all $r,t$ such that $0\le r<t\le T$, and
this equality is satisfied as an equality in $Z_T^{-1}$.
\end{definition}

Similarly to Definition \ref{FWdefsolltion}, the integral with $ds$
in (\ref{intur}) is defined as an element of $Z_T^{-1}$. Definition
\ref{defsolltion} requires for (\ref{parab1}) that this integral
must be equal  to an element of $Z_T^{0}$ in the sense of equality
in $Z_T^{-1}$.

\section{The main result}
For the case of forward SPDEs, the following   result was obtained.
\begin{theorem}
\label{FWTh1}  Assume that $\t=T$ if $\kappa=0$ and $\t\in [0,T)$ if $\kappa\neq 0$.
Problem (\ref{FWparab1})-(\ref{FWparab2}) has a unique solution $u\in Y^1$
for any $\mu\in H^2$, $h_i\in X^1$, and any $\varphi\in X^0$ such that $\oo\varphi\in \U_\t$, where $\oo\varphi(x,t)\defi
\E \varphi(x,t,\o)$. Moreover,
$$ \|u\|_{Y^1}\le
C\left(\|\mu\|_{H^2}+\|\varphi\|_{X^0}+\|\oo\varphi\|_{\U_\t}+\sum_{i=1}^N\|h_i\|_{X^1}\right), $$
Here $c>0$ depends only on $n,T,D,\t$, $\kappa$, $w$, and
on the coefficients of  equation (\ref{FWparab1}).
\end{theorem}
For the case of backward SPDEs, the following result was obtained.
\begin{theorem}
\label{ThB}  
Assume that \baaa \oo\b_i\equiv 0,\quad  \b_i(x,t,\o)\equiv
\oo \b_i(t,\o),\quad  i=1,...,N.
\eaaa 
Let $\e>0$ be given. Assume that $\t=0$ if $\kappa=0$ and $\t\in (0,T]$ if $\kappa\neq 0$.  Then
problem (\ref{parab1})-(\ref{parab2}) has a unique solution
$(u,\chi_1,...,\chi_N)$  in the class $Y^1\times (X^0)^N$,
 for any $\psi\in Z_{T-\e}^2$ and any $\varphi\in X^{0}$ such that $\|\varphi(\cdot,t)\|_{Z_T^0}\equiv 0$ for $t\in[0,\t]$ and
  $\varphi(\cdot,t)\in Z_{T-\e}^0$. In addition, \baaa
\label{3.5} \| u \|_{Y^1}+\sum_{i=1}^N\|\chi_i\|_{X^0}\le C \left(\|
\varphi \| _{X^{0}}+\|\psi\|_{ \Z_{T-\e}^2}\right) \eaaa
for these $\psi$ and $\varphi$.
Here $C>0$ depends only on $\d,T,n,N,D,\e,\t,\lambda$, $\kappa$, $w$, and
on  the
supremums of the coefficients and derivatives  coefficients of  equation (\ref{parab1}).
\end{theorem}
\section{Proofs}
\subsection{Proof of Theorem \ref{FWTh1}}
 We consider the following boundary value problem in $Q$
\baa 
&&d_tu=(\A u+ \varphi)\,dt +\sum_{i=1}^N
[B_iu+h_i]\,dw_i(t), \quad t\ge 0,
\label{F1}
\\
&& u(x,t,\o)\,|_{x\in \p D}=0\label{F2}
\\ &&u(x,0)=\xi(x).\label{F3}
\eaa
 \par
\begin{lemma}
\label{lemma1} Assume that Conditions \ref{cond1}--\ref{condB}
are satisfied.  Then problem   (\ref{F1})-(\ref{F3}) has an unique solution
$u$ in the class $Y^k$ for any $\varphi\in X^{k-2}$,
$h_i\in X^{k-1}$, $\xi\in Z_s^{k-1}$, and \baaa \label{4.2} \| u
\|_{Y^k}\le C \left(\| \varphi \|
_{X^{k-2}}+\|\xi\|_{Z^{k-1}}+ \sum_{i=1}^N\|h_i
\|_{X^{k-1}}\right),\eaaa
where $C > 0$ does not depend on $\xi$ and $\varphi$.
\end{lemma}
\par Note that, in the notations of this lemma, the solution $u=u(\cdot,t)$
is continuous in $t$ in $L_2(\O,\F,\P,Z_T^{k-1})$, since
$Y^k=X^{k}\!\cap \CC^{k-1}$.

The statement of Lemma \ref{lemma1} for $k=1$ can be found in \citet{Roz}, Chapter 3, Section 4.1; the
result represents an analog of the so-called "the first
energy inequality", or "the first fundamental inequality" known for
deterministic parabolic equations (see, e.g., inequality (3.14) from
\citet{la} Ladyzhenskaya (1985), Chapter III).
The statement of Lemma \ref{lemma1} for $k=1$ can be found in \citet{D05}; the result represents an analog of the so-called "the first
energy inequality", or "the first fundamental inequality" known for
deterministic parabolic equations (see, e.g., inequality (3.14) from
\citet{LSU}, Chapter III. It was shown in \citet{D05,Roz} that  $C$ in  Lemma \ref{lemma1}  depends on the
 set of parameters  $$
 n,N,\,\, D,\,\, T,\,\, \,\,\delta, \ \,\,\,\,\\
 \esssup_{x,t,\o,i}\Bigl[| a(x)|+\left|\frac{\p a}{\p x_i}(x)\right|+
|{a_0(x)}|+|\b_i(x,t,\o)|+|\oo\b_i(x,t,\o)|\Bigr].
$$

Let us introduce  operators $\L :H^k\to Y^{k+1}$, $k=0,1$, $L: \W^k\to Y^{k+2}$, $k=-1,0$, and $\H_i:X^k\to Y^{k+1}$, $k=0,1$, such that $\L\xi+L\varphi+\sum _{i=1}^N \H_i h_i=u$,
where $u$ is the solution in $\V$ of  problem (\ref{F1})-(\ref{F3}). By Lemma \ref{lemma1}, all these operators are continuous.

Let us introduce  operators $\oo \L :H^k\to \V^{k+1}$, $k=0,1$, and $\oo L: \W^k\to \V^{k+2}$, $k=-1,0$,  
defined similarly to the   operators $\L :H^k\to \V^{k+1}$, $k=0,1$, and $L: \W^k\to \V^{k+2}$, $k=-1,0$, under the assumptions that
$B_i\equiv 0$, $h_i\equiv 0$, $i=1,...,N$. These linear operators are continuous; see also
e.g. Theorems III.4.1 and IV.9.1 in \citet{LSU} \index{Ladyzhenskaja {\it et al} (1968)} or Theorem III.3.2 in \citet{la}.\index{ Ladyzhenskaya  (1985).}

Let a linear operator $M_0: H^0\to H^1$ be defined
such that $(M_0 \xi)(x)=\int_0^T \WW(t)\oo u(x,t)dt+\kappa\oo u(x,T)$, where $\oo u=\oo \L\xi\in\V^1$.

Further, let linear operator $M: \W^0\to H^1$ be defined
such that $(M \varphi)(x)=\int_0^T \WW(t)\oo u(x,t)dt+\kappa\oo u(x,T)$, where $\oo u=\oo L\varphi\in\V_1$, $\varphi\in \W^0$ (i.e.  $\varphi$ is non-random).

In these notations, $\mu=M_0\oo u(\cdot,0)+M\varphi$ for  a solution $\oo u$ of problem  (\ref{F1})-(\ref{F3})   with $B_i\equiv 0$, $h_i\equiv 0$, $i=1,...,N$, and with nonrandom $\varphi\in \W^0$.

\begin{lemma}
\label{lemma0}\citet{D16}  The linear operator $M_0: H^0\to H^2$  is a continuous bijection; in particular, the inverse operator $M_0^{-1}:H^2\to H^0$ is also continuous. In addition,
the linear operator $M: \U_\t\to H^2$  is  continuous.
\end{lemma}
\par
It can be noted that the classical results for parabolic equations imply that
the operators $M_0: H^k\to H^{k+1}$,  $k=0,1$, and  $M: \W^0\to H^{2}$,  are  continuous for $\kappa=0$, and
the operators $M_0: H^k\to H^{k}$,  $k=0,1$, and  $M: \W^0\to H^{1}$,  are  continuous for $\kappa > 0$; see Theorems III.4.1 and IV.9.1 in \cite{LSU}
or Theorem III.3.2 in \citet{la}. The proof of  continuity of the operator   $M_0: H^0\to H^2$ claimed in Lemma \ref{lemma0} can be found in \citet{D16}. In addition, it was shown in  \citet{D16} that
  the operator $M:\U_\t\to H^2$ is continuous.

Since   problem (\ref{F1})-(\ref{F3}) is linear and the functions $a$ and $a_0$ are non-random, it follows that if  $u=\L\xi+L\varphi+\sum _{i=1}^N \H_i h_i$
and  $\oo u(c,t)\defi \E u(x,t)$, then   $\oo u=\oo\L\xi+\oo L\oo\varphi$, where   $\oo \varphi(x,t)=\E \varphi(x,t)$.
 It follows from the definitions of $M_0$ and $M$ that
\baaa
\mu=M_0\xi+M\oo \varphi.
\eaaa
Hence \baaa
\xi=M_0^{-1}(\mu-M\oo\varphi)\label{xiF}
\eaaa
is uniquely defined, and \baa
u=\L\xi+L\varphi+\sum _{i=1}^N \H_i h_i=\L M_0^{-1}(\mu-M\oo\varphi)+L\varphi+\sum _{i=1}^N \H_i h_i.
\label{usF}
\eaa
 is an unique  solution of problem (\ref{FWparab1})-(\ref{FWparab2}) in $Y^1$.
  By the continuity of the operator $M:\U_\t\to H^2$ and other operators in (\ref{usF}), the
  desired estimate for $u$  follows.
This completes the proof of Theorem \ref{FWTh1}.
$\Box$
\subsection{Proof or Theorem \ref{ThB}}
Let  $\varphi\in X^{-1}$ and $\xi\in Z^0_T$. Consider
the problem \baa
d_tu+\left( \A u+ \varphi\right)dt +
\sum_{i=1}^NB_i\chi_i(t)dt=\sum_{i=1}^N\chi_i(t)dw_i(t), \quad t\le T, \label{bw1}\\
u(x,t,\o)|_{x\in \p D},  \label{bw2}\\
 u(x,T,\o)=\xi(x,\o).  \label{bw3} \eaa
 \par
\begin{lemma}
\label{lemma1B} For $k=1,2$, problem (\ref{bw1})-(\ref{bw3}) has  a unique solution
$(u,\chi_1,...,\chi_N)$ in the class $Y^k\times(X^0)^N$ for any
$\varphi\in X^{k-2}$, $\xi\in Z_T^{k-1}$, and \baaa \label{4.2B} \| u
\|_{Y^k}+\sum_{i=1}^N\|\chi_i\|_{X^{k-1}}\le C \left(\| \varphi \|
_{X^{k-2}}+\|\xi\|_{Z^{k-1}_T}\right), \eaaa where $C>0$ does
not depend on $\varphi$ and $\psi$;
 it depends on $\d,T,n,N,D,$ and on the
supremums of the coefficients and derivatives listed in Condition \ref{cond1}.
\end{lemma}

Note that the solution $u=u(\cdot,t)$
is continuous in $t$ in $L_2(\O,\F,\P,H^{k-1})$, since
$Y^k=X^{1}\!\cap \CC^{0}$.

For $k=1$, the result of Lemma \ref{lemma1B}  can be found in  \citet{D92}  or in \citet{D10} (Theorem 4.2);\index{ Dokuchaev (1992) or Theorem 4.2 from Dokuchaev (2010)).} this is   an analog of the so-called "the first
energy inequality", or "the first fundamental inequality" known for
deterministic parabolic equations; see, e.g., inequality (3.14) in \citet{la}, Chapter III).\index{ from
Ladyzhenskaya (1985), Chapter III).}
\par
For $k=2$, the lemma above represents a reformulation of Theorem 3.1 from \citet{DT}
or Theorem 3.4 in \citet{D10}, or Theorem 4.3 in \citet{D12}; this is  an analog of the so-called "the
second energy inequality", or "the second fundamental inequality"
known for the deterministic parabolic equations (see, e.g.,
inequality (4.56) in \citet{la}, Chapter III. \index{from Ladyzhenskaya (1985), Chapter III).}
 \index{Du
and Tang (2012), or Theorem 3.4 from Dokuchaev (2010) or Theorem 4.3
from Dokuchaev (2012).} In the cited papers, this result was obtained
under some strengthened version of Condition \ref{condA}; this
was restrictive.  In \citet{DT}, this result was obtained
without this restriction, i.e. under Condition \ref{condA}
only.
\begin{remark}\label{remDu} {  Thanks  to Theorem
3.1 from \citet{DT}, Condition 3.5 from  \citet{D11}
and Condition 4.1 from  \citet{D12}   can be replaced by less
restrictive Condition \ref{condA}; all results in \citet{D11,D12}  are still valid.}
\end{remark}
\par
Introduce  operators $L_B:X^{-1}\to Y^1$,
$\L_B:Z^0_s\to Y^1$,  $H_i: X^0\to Y^1$ $\H_i: Z^0_T\to Y^1$ such that $u=L_B\varphi+\L_B\xi$ and $\chi_i=H_i\varphi+\H_i\xi$,  where
 $(u,\chi_1,...,\chi_N)$ is the solution of  problem  (\ref{bw1})-(\ref{bw3})  in the class
$Y^2\times(X^1)^N$. By Lemma \ref{lemma1B}, these linear operators
are continuous.

Let a linear operator $\M_0: H^0\to Z_T^1$ be defined
such that $(\M_0 \xi)(x)=\kappa u(x,0)+\int_0^T \WW(t)\oo u(x,t)dt$, where $u=\L_B\xi\in Y^1$.

Further, let a linear operator $\M: X^0\to  Z_T^1$ be defined
such that $(\M \varphi)(x)=\kappa u(x,0)+\int_0^T \WW(t)\oo u(x,t)dt$, where $u= L_B\varphi\in Y_1$, $\varphi\in X^0$.

In this notations, $\psi=\M_0 \xi+\M\varphi$ for  a solution $u$ of problem  (\ref{bw1})-(\ref{bw3}), i.e. with $\xi=u(\cdot,T)$.

\begin{lemma}
\label{lemmaM0}  The operator $\M_0^{-1}:Z^2_{T-\e}\to Z^0_{T-\e}$ is continuous.
\end{lemma}
\par
{\em Proof of Lemma \ref{lemmaM0}}.
It is known that there exists
an orthogonal basis  $\{v_k\}_{k=1}^\infty$  in $H^0$, i.e. such that
 $$(v_k,v_m)_{H^0}=0,\quad  k\neq m,\quad \|v_k\|_{H^0}=1,$$
and such that $v_k$ is the solution in $ H^1$ of the boundary value problem
\baa
 Av_k=-\lambda_k v_k,\quad  v_k|_{\p D}=0,\label{EP}\eaa
for some $\lambda_k\in\R$, $\lambda_k\to +\infty$ as $k\to +\infty$; see e.g.  Ladyzhenskaya (1985), Chapter 3.4. In other words,
$\lambda_k$ and $v_k$ are the eigenvalues and the corresponding
 eigenfunctions  of the eigenvalue problem (\ref{EP}).

\def\Chi{\Upsilon}
Let $\xi$ and $\psi$ be expanded  as
\baaa
\xi(x,\o)=\sum_{k=1}^\infty \a_k(\o) v_k(x),\quad \psi(x,\o)=\sum_{k=1}^\infty \g_k (\o)v_k(x),\quad \chi_i(x,t,\o)=\sum_{k=1}^\infty \Chi_{ik}(t,\o) v_k(x),
\eaaa
 where
$\{\a_k\}_{k=1}^\infty$ and $\{\g_k\}_{k=1}^\infty$ and square-summable random sequences such that
\baaa
\E \sum_{k=1}^N\a_k^2=\E\|\xi\|_{H^0}^2=\|\xi\|_{\Z_T^0}^2<+\infty,\quad \E \sum_{k=1}^N\g_k^2=\E\|\psi\|_{H^0}^2=\|\psi\|_{\Z_T^0}^2<+\infty.
\eaaa
 By the choice of $\xi$, we have that  $u=\L \xi$. Applying the Fourier method, we obtain that
\baa
u(x,t)=\sum_{k=1}^{\infty}  y_k(t) v_k(x),
\label{sol}\eaa
where $y_k(t)$ are the solutions of backward stochastic differential equations
\baaa
dy_k(t)+\left[-\lambda_k y_k(t)+\sum_{i=1}^N \oo\b_i(t)\Chi_{ik} (t)\right]dt=\sum_{k=1}^N\Chi_{ik}(t)dw_i(t),
\quad y_k(T)=\a_k.
\eaaa
Let \baaa
F(t,s)=\exp\left(-\frac{1}{2}\sum_{i=1}^N\int_s^t\oo\b_i(r)^2dr
+\sum_{i=1}^N\int_s^t\oo\b_i(r)dw_i(r)\right),\quad
\Phi_k(t,s)=e^{-\lambda_k (t-s)}F(t,s).
\eaaa
We have that
\baaa
y_k(t)=\Phi_k(t,0)^{-1}\E\{ \Phi_k(T,0)\a_k|\F_t\}=\E\{ \Phi_k(T,t)\a_k|\F_t\};
\eaaa
see e.g. Proposition 6.2.1 \citet{pham}, p.142.
Hence
\baaa
y_k(t)=e^{-\lambda_k(T-t)}Y_k(t),
\eaaa
where
\baaa
Y_k(t)= \E\left\{F(T,t)\a_k\bigl|\F_t\right\}.
\eaaa
\par
On the other hand,
\baaa
&&
\psi(x)=\sum_{k=1}^{\infty} \g_k v_k(x)= \int_0^T\WW(t)u(x,t)dt+\kappa u(x,0)\\&&=
\sum_{k=1}^{\infty} \int_0^T\WW(t) y_k(t)  v_k(x)dt + \kappa\sum_{k=1}^{\infty}  y_k(0) v_k(x).
\eaaa
Hence $\a_k$ must  be such that $y_k$ satisfy the    equation
 \baa &&\g_k =\int_0^T\WW(t)y_k(t)dt+\kappa y_k(0)
\eaa
which can be rewritten as
\baa &&\g_k=\int_0^T\WW(t)y_k(t)dt+\kappa e^{-\lambda_k T}Y_k(0).
\label{eq_a}\eaa
or
\baa &&\g_k=\int_0^T\WW(t)e^{-\lambda_k (T-t)} Y_k(t)dt+\kappa e^{-\lambda_k T}Y_k(0).
\label{eqaa}\eaa

We consider this as an equation for unknown $\a_k =Y_k(T)=y_k(T)$.
Let us prove  its solvability and uniqueness.
We need to show existence of  $\oo Y_k\in\R$ and $\F_t$-adapted square integrable random processes $\w Y_{kd}$  such that
 \baaa
 \a_k=\oo Y_k+\sum_{d=1}^N \int_0^T \w Y_{kd}(s)dw_d(s).
\eaaa
In this case,   \baaa
Y_k(t)=\oo Y_k+\sum_{d=1}^N \int_0^t \w Y_{kd}(s)dw_d(s).
\eaaa
\par
Let  $\w\g_{kd}$  be  $\F_t$-adapted square integrable random processes such that
\baaa  \g_k=\E \g_k+\sum_{d=1}^N \int_0^{T-\e} \w \g_{kd}(s)dw_d(s),
\eaaa

We assume,
without a loss of generality, that $\lambda_1\ge \lambda$. This is possible because one can always replace the operator $A$
by the operator $A+r I$, for an arbitrarily large $r>0$; in this case, the solution $u$ of the original problem
will be transformed as $u(x,t,\o)e^{-r(t-t)}$, with some adjustment of  $\WW$ and $\varphi$.

Clearly, it suffices to show existence of  $\w Y_{kd}$ and $\oo Y_k\in\R$.
Let $q_k(s)\defi \int_s^T\WW(t)e^{-\lambda_k (T-t)}dt$.
By Condition \ref{cond1}, there exist $\kappa>0$ and $T_1\in [0,T)$ such that
 \baaa
q_k(s)\ge \kappa \int_{T_1\lor s}^Te^{-\lambda_k (T-t)}dt=  \kappa \frac{1-e^{-\lambda_k (T-T_1\lor s)}}{\lambda_k}\ge \kappa \frac{1-e^{-\lambda (T- s)}}{\lambda_k}.
\eaaa
  By the definitions,
\baaa
\E \g_k= \int_0^T\WW(t)e^{-\lambda_k (T- t)}\oo Y_k dt+\kappa  e^{-\lambda_k T} \oo Y_k,
\eaaa
and \baaa \oo Y_k=  \left[q_k(0)+\kappa  e^{-\lambda_k T}\right]^{-1}\E\g_k.
\eaaa
Hence
\baaa
\oo Y_k^2\le \const \lambda_k^2\E \g_k^2.
\eaaa

In addition,
\baaa \int_0^{T-\e}\w \g_{kd}(s)dw_d(s)=\int_0^T\WW(t)e^{-\lambda_k (T-t)}dt\int_0^t \w Y_{kd}(s)dw_d(s)\\
=\int_0^T\w Y_{kd}(s)dw_d(s) \int_s^T\WW(t)e^{-\lambda_k (T-t)}dt .
\eaaa
Hence
\baaa \w Y_{kd}(s)=\Ind_{\{s\le T-\e\}}\frac{\w \g_{kd}(s)}{q_k(s)}
\eaaa
and
\baaa
|\w Y_{kd}(s)|\le  \Ind_{\{s\le T-\e\}}\frac{\lambda_k|\w \g_{kd}(s)|}{\kappa(1-e^{-\lambda (T- s)})}.
\eaaa
 Hence equation (\ref{eqaa}) has an unique solution $\a_k\in L_2(\O,\F_{T-\e},\O)$ and  there
exists  $C>0$ such that
\baa
\E\sum_{k=1}^{\infty} \a_k^2\le C\E \sum_{k=1}^{\infty} \g_k^2 \lambda_k^2.
\label{3}\eaa

Further, we have that \baaa
A\psi=\sum_{k=1}^{\infty} \g_k A v_k =-\sum_{k=1}^{\infty} \g_k \lambda_k v_k\quad\hbox{a.s.}
\eaaa
and
\baaa
\E\|A\psi\|_{H^0}^2=\E\sum_{k=1}^{\infty} \g_k^2 \lambda_k^2.
\label{Amuf}
\eaaa
Hence (\ref{3}) can be rewritten as
\baa
\E\|\xi\|_{Z_{T-\e}^0}^2\le C_1\E\|A\psi\|_{H^0}^2\le C \| \psi\|_{Z^2_{T-\e}}^2
\label{3n}
\eaa
for some $C_1>0$ and $C>0$ that are
independent on $\psi$. Thus, (\ref{3n})   implies that the operator $\M_0^{-1}:Z^2_{T-\e}\to Z^0_{T-\e}$ is continuous.
This completes the proof of Lemma \ref{lemmaM0}. $\Box$

\begin{lemma}
\label{lemma2B} Solution of problem (\ref{bw1})-(\ref{bw3}) for
$\varphi=0$ and $\xi\in Z_T^0$ is such that   \baaa \label{4nB} \| u
(\cdot,0)\|_{Z_0^2}\le C\|\xi\|_{Z^{0}_T},\eaaa where $C>0$ does
not depend on $\xi$;
 it depends on $\d,T,n,N,D,$ and on the
supremums of the coefficients and derivatives listed in Condition \ref{cond1}.
\end{lemma}

 {\em Proof of Lemma \ref{lemma2B}}.
 We have that $u(\cdot,0)=\sum_{k=1}^\infty v_k(x)y_k(0)$. Hence
 \baaa
 \|Au(\cdot,0)\|_{Z_T^0}=\E\sum_{k=1}^\infty \lambda_k^2 y_k(0)^2=\sum_{k=1}^\infty \lambda_k^2 e^{-\lambda_k T}[\E (F(T,0)\a_k)]^2\\
 \le \sum_{k=1}^\infty \lambda_k^2e^{-\lambda_k T}[\E (F(T,0)\a_k)]^2
 \eaaa
We have that $[\E (F(T,0)\a_k)]^2\le c_F \E\a_k^2$, where $c_F\defi \E F(T,0)^2$.
Hence
\baaa
 \|Au(\cdot,0)\|_{Z_T^0}
 \le c_F\sum_{k=1}^\infty \lambda_k^2e^{-\lambda_k T}\E \a_k^2.
 \eaaa
Hence
\baaa
 \|Au(\cdot,0)\|_{Z_T^0}
 \le C\|\xi\|_{Z_T^0},
 \eaaa
where $C>0$   does
not depend on $\xi$;
 it depends on $\d,T,n,N,D,$ and on the
supremums of the coefficients and derivatives listed in Condition \ref{cond1}.
Further, for any $\lambda\in \R$,  we have that $h\defi A u(\cdot,0)+ \lambda u(\cdot,0)\in H^0$. By the properties of the elliptic equations,
it follows  that there exists $\lambda\in\R$ and $c=c(\lambda)>0$ such that \baa
\|u(\cdot,0\|_{H^2}\le c\|h\|_{H^0}\le c(\|A u(\cdot,0)\|_{H^0}+\|\lambda u(\cdot,0)\|_{H^0})\quad\hbox{a.s.};
\label{2ndfund}\eaa
see e.g. Theorem II.7.2 and Remark II.7.1 in Ladyzhenskaya (1975), or  Theorem III.9.2 and Theorem  III.10.1 in Ladyzhenskaya and Ural'ceva  (1968).
By (\ref{2ndfund}), we have that  \baa
\|u(\cdot,0)\|_{Z^0_2}\le  c_1(\|A u(\cdot,0)\|_{Z_T^0}+\|\xi\|_{Z_T^0})\le c_2\|\xi\|_{Z^0_T}
\label{2ndfund2}\eaa
for some $c_i>0$ that are independent on $\xi$ but depend on the parameters of the SPDE. This completes the proof of Lemma \ref{lemma2B}. $\Box$

\begin{lemma}
\label{lemmaM00}  The operators $\M_0:Z^2_{T}\to Z^0_{T}$ and $\M_0:Z^2_{T-\e}\to Z^0_{T-\e}$ are 
 continuous.
\end{lemma}
\par
{\em Proof of Lemma \ref{lemmaM00}}.
 By Lemmas  \ref{lemma1B}-\ref{lemma2B}, the operator $\M_0:Z^0_{T}\to Z^2_{T}$ is continuous.
 Further, let $\xi\in Z_{T-\e}^0$.
Since $\chi_i(t)=0$ for $t>T-\e$ and $\chi_i=\H_i\xi$, $\xi\in  Z^0_{T-\e}$,
 it follows that the operator $\M_0:Z^2_{T-\e}\to Z^0_{T-\e}$ is continuous.
This completes the proof of Lemma \ref{lemmaM00}. $\Box$

Let $\X_\t$ be the space of all $\varphi\in X^0$ such that $\|\varphi(\cdot,t)\|_{Z_0^T}=0$ for $t\in[0,\t)$, with the norm $\|\varphi\|_{\X_\t}=\|\varphi\|_{X^0}$.

Let $\Y_\e$ be the space of all $\varphi\in X^0$ such that
  $\varphi(\cdot,t)\in Z_{T-\e}^0$, with the norm $\|\varphi\|_{\Y_\e^0}=\|\varphi\|_{X^0}$.
  
We consider  $\X_\t \cap \Y_\e$ as a linear normed space with the norm $\|\varphi\|_{\Y_\e^0}=\|\varphi\|_{X^0}$.
Clearly, this is a Banach space.

\begin{lemma}
\label{lemmaM}  The
operators $\M:\X_\t\to \Z_T^2$, $\M:\Y_\e\to \Z_{T-\e}^1$, and $\M:\X_\t\cap \Y_\e  \to \Z_{T-\e}^2$ are continuous.
\end{lemma}
\par
{\em Proof of Lemma \ref{lemmaM}}.
 Let us show that the operator $\M:\X_\t\to Z_T^2$ is continuous.
 By Lemma \ref{lemma1B}, the operator $\M:X^0\to Z_T^2$ is continuous if $\kappa=0$; in this case, we can select $\t=T$ and $X_\t^0=X^0$.

Let us show  that the operator $\M:\X_\t\to Z^2_T$ is continuous for the case where $\kappa\neq 0$ and $\t\neq 0$.  Without a loss of generality, let us assume  that  $\kappa=1$, $\mu =\M \varphi=u(x,0)$, and $\WW(t)\equiv 0$; it suffices  because  the boundary value problem is linear.

Let  $\varphi\in X_\t^{0}$  and $u=\L_B\varphi$.
 For $\w\t\in (0,\t]$, consider a modification of
problem  (\ref{bw1})-(\ref{bw3}) \baaa \label{semi}
\begin{array}{ll}
d_t\w u+\left[A u  +
\sum_{i=1}^NB_i\chi_i(t)\right]dt=\sum_{i=1}^N\chi_i(t)dw_i(t), \quad t\in (0,\w\t),\\
\w u(x,t,\o)|_{x\in \p D}, \\
 \w u(x,\w\t,\o)=u(x,\w\t,\o).
\end{array}
 \eaaa
 By the semi-group property of backward SPDEs from Theorem
6.1 from \cite{D10}, we obtain that $\w u|_{t\in[0,\w\t]}=u|_{t\in[0,\t]}$.
By Lemma
\ref{lemma1B}, we have for $\tau\in[\w\t,\t]$
that
\baaa \inf_{t\in[\w t,\t]} \|
u(\cdot,t)\|^2_{Z_t^2}\le \frac{1}{\t-\w\t}
\int_{\w\t}^\t \|u(\cdot,t)\|_{Z_t^2}^2dt\le
\frac{C_1}{\t-\w\t}\|\varphi\|^2_{X^0},\eaaa
where $C_1>0$ is independent of $\varphi$.
 Hence the exists $\tau\in[\w t,\t]$ such that
 \baaa  \|u(\cdot,\tau)\|^2_{Z_\tau^2}\le
C\|\varphi\|^2_{X^0},\eaaa
 where $C=2C_1(\t-\w\t)^{-1}$. Repeating the proof
 for   the continuity of the operator  $\M_0:Z_T^2\to  Z_T^0$
 to the new time interval $[0,\tau]$, we obtain that
 the operator $\M:\X_\t \to Z_T^2$ is continuous.

 Further, let $\varphi\in \Y_{\e}$. In this case, 
 $\chi_i(t)=0$ for $t>T-\e$ and $\chi_i=H_i\varphi$ and $\M\varphi\in  Z^1_{T-\e}$.
 it follows that the operator $\M:\Y_{\e}\to Z^1_{T-\e}$ is continuous.
 Then the continuity of the operator  $\M:\X_\t\cap \Y_\e  \to Z_{T-\e}^2$ follows.
This completes the proof of Lemma \ref{lemmaM}. $\Box$

\vspace{4mm}

We are now in the position to complete the proof of Theorem \ref{ThB}.
By the assumptions, $\varphi\in \X_\t\cap \Y_\e$.
 It follows from the definitions of $\M_0$ and $\M$ that
\baaa
\psi=\M_0\xi+\M\varphi.
\eaaa
Since the operator $\M:X_\t^0 \to Z_T^2$  and $\M_0^{-1}: Z_T^2 \to Z_T^0$ are  continuous, it follows that   $\M\varphi\in Z_T^2$ and
\baa
\xi=\M_0^{-1}(\psi-\M\varphi)\label{xi}
\eaa
is uniquely defined. Hence \baa
&&u=\L_B\xi+L_B\varphi=\L_B \M_0^{-1}(\psi-\M\varphi)+L_B\varphi,\nonumber\\&&
\chi_i=H_i\varphi+\H_i\xi=H_i\varphi_i+\H_i\M_0^{-1}(\psi-\M\varphi),\quad i=1...,N,
\label{us}
\eaa
 is a unique  solution of problem (\ref{parab1})-(\ref{parab2}) in $Y_1$.
   By the continuity of the operator $\M:\X_\t\cap \Y_\e  \to Z_{T-\e}^2$  and other operators in (\ref{us}), the
  desired estimate for $u$  follows.  This completes the proof of Theorem \ref{ThB}.
$\Box$
\section*{Discussion and future development}
\begin{enumerate}
\item Consider the case where $\kappa=0$ and $\WW\equiv 1$.
By Theorem \ref{ThB}, the corresponding boundary problem is well-posed as a backward SPDE (or BSPDEs), even
without  a Cauchy condition that is usually associated with BSPDEs. This gives  a new information about
the nature of BSPDEs.
\item
The current proof is based on eigenfunctions expansion; we think that it is possible 
to obtain a similar result using a different proof for  the case where  the operator $A$ is not necessarily symmetric and has coefficients depending on time. We leave this for the future research.
\item 
\item  It would be interesting to extend the result for the case of backward SPDEs on the case where $\b_i=\b_i(x,t,\o)$ or where the operators $B_i$ are differential operators. So far, the approach based on eigenfunctions expansion  used here does not allow this. 
\end{enumerate}

\end{document}